\documentclass[12pt]{article}

\usepackage{amsthm,amsmath,amssymb}

\date{}

\iftrue \makeatletter
\@ifundefined{every@math@size}{}{%
\addto@hook\every@math@size{\dch@scr@hook}
\def\dch@scr@adjust{\@ifundefined{dch@sizet\f@size}%
  {\expandafter\dch@set@script\csname dch@sizet\f@size\endcsname}%
  {\csname dch@sizet\f@size\endcsname}}
\def\dch@set@script#1{\begingroup %
  \frozen@everymath{}%
  \let#1\@empty \let\dch@do@one\relax
  \dch@set@one\scriptscriptstyle\scriptscriptfont#1\ssf@size
  \dch@set@one\scriptstyle\scriptfont#1\sf@size
  \dch@set@one\textstyle\textfont#1\f@size
  \endgroup #1} %
\def\dch@set@one#1#2#3#4{%
  \@ifundefined{dch@size#4}%
   {\expandafter\xdef\csname dch@size#4\endcsname{%
      \fontdimen13\the#2\tw@=\the\fontdimen13#2\tw@
      \fontdimen14\the#2\tw@=\the\fontdimen14#2\tw@
      \fontdimen15\the#2\tw@=\the\fontdimen15#2\tw@
      \fontdimen16\the#2\tw@=\the\fontdimen16#2\tw@
      \fontdimen17\the#2\tw@=\the\fontdimen17#2\tw@}%
  }{\csname dch@size#4\endcsname}%
  \setbox\z@\hbox{$#1H_2$}\@tempdima\dp\z@
  \setbox\z@\hbox{$#1H_2^{+\vrule \@height 1em}$}%
   \ifdim\@tempdima<\dp\z@
    \advance\@tempdima\dp\z@ \divide\@tempdima\tw@ %
    \@tempdimb\dp\z@ \advance\@tempdimb-\@tempdima %
    \advance\@tempdimb\ht\z@ \advance\@tempdimb-1em %
    \xdef#3{#3\dch@do@one#2{\the\@tempdimb}{\the\@tempdima}}%
  \fi}
\def\dch@do@one#1#2#3{%
  \fontdimen13#1\tw@#2\relax
  \fontdimen14#1\tw@\fontdimen13#1\tw@ \fontdimen15#1\tw@\fontdimen13#1\tw@
  \fontdimen\sixt@@n#1\tw@#3\fontdimen17#1\tw@\fontdimen\sixt@@n#1\tw@}%
\let\dch@scr@hook\dch@scr@adjust
\ifx\glb@currsize\f@size \dch@scr@adjust \fi}%
\makeatother \fi

\renewcommand{\baselinestretch}{1.08}\normalsize

\newtheorem{thm}{Theorem}

\newtheorem{lemma}{Lemma}

\def\dfrac#1#2{\lower0.15ex\hbox{\large$\frac{#1}{#2}$}}
\def\({\bigl(}
\def\){\bigr)}

\def\B{{\cal B}}
\def\Bn{\B_n}
\def\T{\mathcal{T}}
\def\Tmn{\T_{m,n}}

\def\vol{\operatorname{vol}}
\def\relvol{\operatorname{\nu}}
\def\Reals{{\mathbb{R}}}
\def\Ints{{\mathbb{Z}}}

\def\nicebreak{\vskip 0pt plus 50pt\penalty-300\vskip 0pt plus -50pt }

\begin{document}

\title {The asymptotic volume of the Birkhoff polytope}

\author{
E.~Rodney~Canfield\thanks
 {Research supported by the NSA Mathematical Sciences Program.} \\
\small Department of Computer Science\\[-0.8ex]
\small University of Georgia\\[-0.8ex]
\small Athens, GA 30602, USA\\[-0.3ex]
\small\texttt{ercanfie@uga.edu}
\and
\vrule width0pt height3ex
Brendan~D.~McKay\vrule width0pt height2ex\thanks
 {Research supported by the Australian Research Council.}\\
\small Department of Computer Science\\[-0.8ex]
\small Australian National University\\[-0.8ex]
\small Canberra ACT 0200, Australia\\[-0.3ex]
\small\texttt{bdm@cs.anu.edu.au}
}

\maketitle

\vspace*{-5mm}
\begin{abstract}
Let $m,n\ge 1$ be integers.
Define $\Tmn$ to be the \textit{transportation polytope\/}
consisting of the $m\times n$ non-negative real matrices
whose rows each sum to~$1$ and whose columns each sum to~$m/n$.
The special case $\Bn=\T_{n,n}$ is the much-studied
\textit{Birkhoff-von Neumann polytope\/} of doubly-stochastic matrices.  
Using a recent asymptotic enumeration of non-negative integer matrices
(Canfield and McKay, 2007), we determine the asymptotic volume of
$\Tmn$ as $n\to\infty$ with $m=m(n)$ such that $m/n$ neither decreases
nor increases too quickly.  In particular, we give an asymptotic 
formula for the volume of $\Bn$.
\end{abstract}

\nicebreak
\section{Introduction}\label{s:intro}

Let $m,n\ge 1$ be integers.
Define $\Tmn$ to be the \textit{transportation polytope\/}
consisting of the $m\times n$ non-negative real matrices
whose rows each sum to~$1$ and whose columns each sum to~$m/n$.
The special case $\Bn=\T_{n,n}$ is the famous
\textit{Birkhoff-von Neumann polytope\/} of doubly-stochastic matrices.

It is well known (see Stanley~\cite[Chap. 4]{Stanley} for basic
theory and references) that $\Tmn$ spans an
$(m{-}1)(n{-}1)$-dimensional affine subspace of $\Reals^{m\times n}$.
The vertices of $\Tmn$ were described by Klee and
Witzgall~\cite{Klee} and are moderately complicated.  The special
case of~$\Bn$ is however very simple: the vertices are precisely
the $n\times n$ permutation matrices.

Two types of volume are customarily defined for such
polytopes.  We can illustrate the difference using the example
\[
\B_2 = \biggl\{
    \Bigl(\,\begin{matrix} z&\!\!1{-}z\\1{-}z&\!\!z\end{matrix}\,\Bigr)
       \;\biggm|\; 0\le z\le 1\,\biggr\}
    = \biggl[\,
        \Bigl(\,\begin{matrix} 0&\!1\\1&\!0\end{matrix}\,\Bigr),
        \Bigl(\,\begin{matrix} 1&\!0\\0&\!1\end{matrix}\,\Bigr)
      \,\biggr],
\]
where the last notation indicates a closed line-segment
in~$\Reals^{2\times 2}$.
The length
of this line-segment is the \textit{volume} $\vol(B_2)=2$.
We can also consider the lattice induced by $\Ints^{2\times 2}$ on the
affine span of $\B_2$: this consists of the points
$\(\begin{smallmatrix} z&1{-}z\\1{-}z&z\end{smallmatrix}\)$
for integer~$z$.  The polytope~$\B_2$ consists of a single basic
cell of this lattice, so it has \textit{relative volume\/}
$\relvol(\B_2)=1$.
In general, $\vol(\Tmn)$ is the volume in units of the ordinary
$(m{-}1)(n{-1})$-dimensional Lebesgue  measure,
while $\relvol(\Tmn)$ is the volume in units of basic cells
of the lattice
induced by $\Ints^{m\times n}$ on the affine span of~$\Tmn$.

\bigskip
\begin{lemma}\label{ratio}
 For $m,n\ge 2$,
 $\vol(\Tmn) = m^{(n-1)/2} n^{(m-1)/2}\relvol(\Tmn)$.
\end{lemma}
\begin{proof}
 This is established in~\cite[Theorem~3]{diaconis}.
 Also see the Appendix of~\cite{ChanRobbins}.
\end{proof}

Next, define the function $H_{m,n}:\Ints\to\Ints$ by
\[
  H_{m,n}(z) = \bigl| z\Tmn\cap\Ints^{m\times n} \bigr|.
\]
Clearly $z\Tmn\cap\Ints^{m\times n}$ is the set of $m\times n$ non-negative
integer matrices with row sums equal to $z$ and column sums equal
to~$zm/n$.  This set is non-empty when $zm/n\in\Ints$; that is,
when $z$ is a multiple of $z_0=n/\gcd(m,n)$. The base case
$z=z_0$ corresponds to an expanded polytope $z_0\Tmn$
whose vertices are integral~\cite[Cor.~1]{Klee}.
Therefore, by the celebrated theorem of Ehrhart (see~\cite{Stanley}),
there are constants $c_i(m,n)$ for $i=0,1,\ldots,(m{-}1)(n{-}1)$
such that
\begin{equation}\label{ehrhart}
  H_{m,n}(z) = \begin{cases}\;\displaystyle\sum_{i=0}^{(m{-}1)(n{-}1)}
     c_i(m,n) z^{(m{-}1)(n{-}1)-i}, &\text{if $z_0$ divides $z$;}\\[1ex]
     \;0,&\text{otherwise.}\end{cases}
\end{equation}
This is the \textit{Ehrhart pseudo-polynomial\/} of $\Tmn$.
Applying~\cite[Prop.~4.6.30]{Stanley} to $z_0\Tmn$, we find that
\begin{equation}\label{leading}
  \relvol(\Tmn) = c_0(m,n).
\end{equation}

\medskip

We turn now to asymptotics.  Our main tool will be the following
theorem of the present authors~\cite{CMinteger}.

\begin{thm}\label{Mvalue}
Suppose $m=m(n)$, $s=s(n)$ and $t=t(n)$ are positive integer
functions such that $ms=nt$.  
Let $M(m,s;n,t)$ be the number of $m\times n$ non-negative
integer matrices with row sums equal to~$s$ and column sums
equal to~$t$.
Define $\lambda=\lambda(n)$ by $ms=nt=\lambda mn$.
Let $a,b>0$ be constants
such that $a+b<\frac12$.  Suppose that $n \rightarrow \infty$
and that, for large $n$,
\begin{equation} \label{Hyp}
\frac{(1+2\lambda)^2}{4\lambda(1+\lambda)}
   \left( 1 + \frac{5m}{6n} + \frac{5n}{6m} \right)
\le a \log n.
\end{equation}
Then
\begin{align*}
M(m,s;n,t) &=
  \frac{\displaystyle\binom{n{+}s{-}1}{n-1}^{\!m}
           \binom{m{+}t{-}1}{m-1}^{\!n}}
          {\displaystyle\binom{mn{+}\lambda mn{-}1}{mn-1}}
    \exp\(\, \dfrac12 + O(n^{-b}) \).\quad\qedsymbol
\end{align*}

\end{thm}

\noindent Using this result, we can prove the following theorem
concerning the volumes of $\Tmn$ and~$\Bn$.

\begin{thm}\label{main}
  Let $a,b>0$ be constants such that $a+b<\tfrac12$.
  Then
 \[
  \vol(\Tmn) =
     \frac{1}{(2\pi)^{(m+n-1)/2} n^{(m-1)(n-1)}}
     \exp\Bigl(\dfrac13 + mn - \frac{(m-n)^2}{12mn} + O(n^{-b})\Bigr)
 \]
 when $m,n\to\infty$ in such a way that\/
 $\max\(\dfrac mn,\dfrac nm\)\le\dfrac65\,a\log n$.
 In particular, for any $\epsilon>0$
 \[
   \vol(\Bn) = \frac{1}{(2\pi)^{n-1/2} n^{(n-1)^2}}
       \exp\Bigl(\dfrac13 + n^2 + O(n^{-1/2+\epsilon})\Bigr)
 \]
 as $n\to\infty$.
\end{thm}
\begin{proof}
{}From~\eqref{ehrhart} and~\eqref{leading}, we have
\begin{equation}\label{limit}
  \relvol(\Tmn) = \lim_{z\to\infty} \frac{H_{m,n}(z)}{z^{(m-1)(n-1)}}
      = \lim_{\lambda\to\infty} 
        \frac{M(m,\lambda n; n,\lambda m)}{(\lambda n)^{(m-1)(n-1)}}\,,
\end{equation}
where we restrict $z$ to multiples of $z_0$ and $\lambda$ to
multiples of $z_0/n$.
If $a'>a$ and $a'+b<\tfrac12$, then the left side of~\eqref{Hyp}
is less than $a'\log n$ for sufficiently large~$\lambda$.
Thus the conditions for Theorem~\ref{Mvalue} hold.  It remains
to apply that theorem to~\eqref{limit} using Stirling's formula,
and to infer the value of $\vol(\Tmn)$ using~Lemma~\ref{ratio}.
\end{proof}

It is of interest to note that the same asymptotic formula for the
volume (except for the error term) follows from the
estimate of $M(m,s;n,t)$ that Diaconis and Efron
proposed without proof in 1985~\cite{diaconis}.

\begin{table}[ht]
\centering
\begin{tabular}{c|c}
 $n$ & estimate/actual \\
 \hline
 1 & 1.51345 \\
 2 & 1.20951 \\
 3 & 1.25408 \\
 4 & 1.22556 \\
 5 & 1.19608 \\
 6 & 1.17258 \\
 7 & 1.15403 \\
 8 & 1.13910 \\ 
 9 & 1.12684 \\
 10 & 1.11627
\end{tabular}
\caption{Accuracy of Theorem~\ref{main} for $\vol(\Bn)$.\label{tab}}
\end{table}

Exact values of $\vol(\Bn)$ are known up to $n=10$~\cite{BP2003}.
In Table~\ref{tab} we compare the exact values to the approximation
given in Theorem~\ref{main}.  It appears that the true magnitude
of the error term might be~$O(n^{-1})$.
This would indeed be the case if the well-tested conjecture
made in~\cite{CMinteger} about the value of $M(n,s;n,t)$ was true.
The same conjecture implies a value of
$\vol(\Tmn)$ with relative error $O\((m+n)^{-1}\)$ for all~$m,n$.

Recently, a summation with $O\(n^n n!\)$ terms was
found for $\vol(\Bn)$~\cite{deloera}.  Whether it is useful for asymptotics
remains to be seen.

\nicebreak
\renewcommand{\baselinestretch}{1.0}

\end{document}